\documentclass[leqno,11pt]{amsart}
\usepackage{amsmath,amssymb}
\usepackage{tikz, tkz-fct}
\usepackage[margin=1.25in]{geometry}
\usepackage[colorlinks,pagebackref,hypertexnames=false]{hyperref}
\newtheorem{theorem}{Theorem}
\theoremstyle{definition}
\newtheorem{example}[theorem]{Example}

\title{Examples of equivariant Lagrangian mean curvature flow}

\author{Jason D. Lotay}
\address{University of Oxford, U.K.}
\urladdr{\href{http://people.maths.ox.ac.uk/lotay/}{http://people.maths.ox.ac.uk/lotay/}}
\email{jason.lotay@maths.ox.ac.uk}
 
\begin{document}

\begin{abstract}{} \noindent In this expository note we describe important examples of Lagrangian mean curvature flow in $\mathbb{C}^2$ which are invariant under a circle action.  Through these examples, we see compact and non-compact situations,  long-time existence, singularities forming via explicit models, and significant objects in Riemannian and symplectic geometry, including the Clifford torus, Chekanov torus, Whitney sphere and Lawlor necks.
\end{abstract}

\maketitle

\section{Introduction}

Lagrangian mean curvature flow has generated major interest from several viewpoints, due to its connections with Riemannian geometry (particularly calibrated geometry), symplectic topology, gauge theory, Calabi--Yau (and, more generally, K\"ahler--Einstein) manifolds and Mirror Symmetry.  In particular, Lagrangian mean curvature flow has the potential to lead to striking applications in diverse areas.  

However, there are relatively few cases in which the Lagrangian mean curvature flow is explicitly understood.  An exception is the setting of equivariant Lagrangian mean curvature flow in $\mathbb{C}^n$, which was first studied in \cite{Grohetal,NevesMonotone,NevesSingularities}.    In particular, understanding  equivariant flows in $\mathbb{C}^2$ was crucial in the ground-breaking result by Neves \cite{NevesSingularities} that \emph{any} embedded Lagrangian in a Calabi--Yau 2-fold is Hamiltonian isotopic to a Lagrangian which develops a finite-time singularity under Lagrangian mean curvature flow.   

In the case of $\mathbb{C}^2$ (with its standard symplectic form), which will be the focus of this article, circle-invariant Lagrangian surfaces $L\subseteq\mathbb{C}^2$ are given by curves $\gamma:I\subseteq\mathbb{R}\to\mathbb{C}$ as follows:
\begin{equation}\label{eq:L.equiv}
L=\{(\gamma(s)\cos \phi,\gamma(s)\sin \phi)\in\mathbb{C}^2:s\in I,\phi\in[0,2\pi)\}.
\end{equation}
This is the heart of the advantage of studying equivariant Lagrangian mean curvature flow, in that it reduces to a flow of curves in $\mathbb{C}$.  (It is worth emphasising that points on $\gamma$ define circles in the Lagrangian $L$ in \eqref{eq:L.equiv}, except the origin which gives just a point in $L$.)

  Recent exciting progress has been made in obtaining a detailed understanding of equivariant Lagrangian mean curvature flow in $\mathbb{C}^2$ (and, more generally, in $\mathbb{C}^n$), namely in  \cite{ELS,SmoczykWhitney,Su,Viana,Wood}.  We shall briefly describe the main outcomes of this work in this article, which we organize into three cases  that depend on the curve $\gamma$ in $\mathbb{C}$ defining our equivariant Lagrangian:
\begin{itemize}
\item Lagrangian tori $T^2$, given by embedded closed curves;
\item Lagrangian spheres $S^2$, given by immersed closed curves;
\item Lagrangian cylinders $S^1\times\mathbb{R}$, given by embedded open curves.
\end{itemize}
(One may also have planes $\mathbb{R}^2$, given by certain open curves, but these have not yet been studied.)  

  The examples described here provide important results for Lagrangian mean curvature flow, and it certainly motivates future research in the equivariant setting to obtain further progress in our understanding of the flow.

\section{Preliminaries}

We begin with some fundamental notions we require  to describe the examples.

If $z_1=x_1+iy_1$, $z_2=x_2+iy_2$ are standard coordinates on $\mathbb{C}^2$, we have the standard symplectic form given by
\begin{equation*}
\omega=\frac{i}{2}(d z_1\wedge d\bar{z}_1+d z_2\wedge d\bar{z}_2)=d x_1\wedge d y_1+d x_2\wedge d y_2.
\end{equation*}
One may then check that surfaces $L$ as in \eqref{eq:L.equiv} are Lagrangian, namely $\omega|_L\equiv 0$.  Lagrangian mean curvature flow in $\mathbb{C}^2$ (and any K\"ahler--Einstein manifold) is  the   mean curvature flow with a Lagrangian initial condition:
\begin{equation}\label{eq:LMCF}
\frac{\partial L}{\partial t}=H,
\end{equation}
where $H$ is the mean curvature vector of $L$.

If we take a circle-invariant Lagrangian $L$ as in \eqref{eq:L.equiv}, then Lagrangian mean curvature flow $L_t$ preserves the circle-invariance, so we may write the flow \eqref{eq:LMCF} as a flow of curves $\gamma_t$ in $\mathbb{C}$:
\begin{equation}\label{eq:curve.flow}
\frac{\partial\gamma}{\partial t}=\kappa-\frac{\gamma^{\perp}}{|\gamma|^2},
\end{equation}
where $\kappa$ is the curvature vector of $\gamma$ with respect to the Euclidean metric on $\mathbb{C}$, and $\gamma^{\perp}$ is the projection of  $\gamma$ to the normal direction to $\gamma$.  

The key to understanding any geometric flow is to study the formation of singularities.  There are two types of singularities in Lagrangian mean curvature flow at a finite time $t=T$, determined by the behaviour of the second fundamental form $A_t$ of the flow $L_t$, which must blow-up at the singularity.  
\begin{itemize}
\item Type I: singularities where $\lim_{t\nearrow T} \sup_{L_t}|A_t|^2(T-t)<\infty$.
\item Type II: any other singularity at $t=T$.
\end{itemize}
The simplest examples of Type I singularities are given by \emph{self-shrinkers}: these are solutions to \eqref{eq:LMCF} that simply shrink by dilations under the flow.  

There are also two main ways of studying these singularities via blow-up.
\begin{itemize}
\item Type I blow-up.  For a positive sequence $\sigma_i\to\infty$ and $w\in\mathbb{C}^2$, we define 
\begin{equation*}
L^i_s=\sigma_i(L_{T+\sigma_i^{-2}s}-w)\quad \forall\, s\in [-\sigma_i^2T,0).
\end{equation*}   
The sequence $L^i_s$ subconverges weakly (i.e.~as a Brakke flow) as $i\to\infty$, to a limit flow $L^{\infty}_s$ in $\mathbb{C}^2$ for all $s<0$, which is called a Type I blow-up (or tangent flow) at $(w,T)$.  A Type I blow-up is a self-shrinker and it provides a ``first approximation'' to the nature of the singularity. 
\item Type II blow-up.  Suppose we have a sequence $(w_i,t_i)\in \mathbb{C}^2\times (0,T)$ such that $t_i\to T$ and
 $\sigma_i:=A_{t_i}(w_i)=\sup\{|A_t(z)|\,:\,z\in L_t,\,t\leq t_i\}>0$.
If  
\begin{equation*}
L^i_s=\sigma_i(L_{t_i+\sigma_i^{-2}s}-w_i) \quad\forall\, s\in[-\sigma_i^2t_i,0),
\end{equation*}
then the sequence $L^i_s$ subconverges as $i\to\infty$ and it will define a smooth  solution $L^{\infty}_s$ in $\mathbb{C}^2$ for all $s<0$, which we call a  Type II blow-up. Type II blow-ups give more refined information about the singularity formation.
\end{itemize}
\noindent Since Type I and Type II blow-ups are flows defined for all negative times, they are examples of ancient solutions to the flow.

 Suppose $L_t$ is a Type II blow-up.  For a sequence $\lambda_i\to\infty$, we let
\begin{equation*}
L_s^i :=\lambda_i^{-1} L_{\lambda_i^2s} \qquad  \forall s<0,
\end{equation*}
and define the \emph{blow-down} $L^{\infty}_s$ as a subsequential limit of the sequence $L^i_s$.  A blow-down of a Type II blow-up is a self-shrinker for Lagrangian mean curvature flow in $\mathbb{C}^2$ and will   ``approximate'' the Type II blow-up.  Recently, classification results have been obtained for Type II blow-ups in Lagrangian mean curvature flow in terms of their blow-downs \cite{LLS}.

\section{Tori: embedded closed curves}

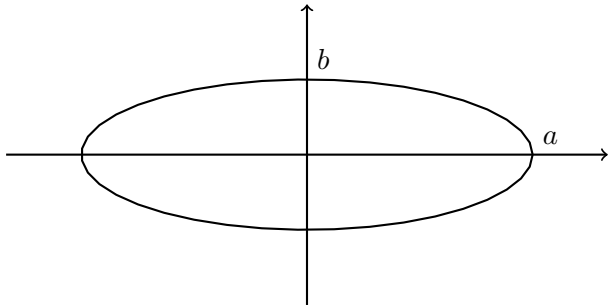
\begin{figure}[h]
\begin{center}
\begin{tikzpicture}[xscale=1,yscale=1]
\draw[thick,->] (0,-2) -- (0,2);
\draw[thick,->] (-4,0) -- (4,0);
\node[above right] at (3,0) {$a$};
\node[above right] at (0,1) {$b$};
\draw[thick,domain=0:2*pi,samples=40] plot ({3*cos(\x r)},{1*sin(\x r)}) ;
\end{tikzpicture}
\caption{An ellipse defining an embedded Lagrangian 2-torus $T^2$}\label{fig:ellipse}
\end{center}
\end{figure}

In this section we examine embedded, equivariant, Lagrangian 2-tori $L$ in $\mathbb{C}^2$ which are defined by embedded closed curves $\gamma$ in $\mathbb{C}$ that enclose the origin as in Figure \ref{fig:ellipse}.

\begin{example}[Clifford torus/circle] The case when $\gamma$ is a circle, i.e.~when $a=b$ in Figure \ref{fig:ellipse}, corresponds to $L$ being the well-known Clifford torus (which is Lagrangian) in $\mathbb{C}^2$.  The Clifford torus is a self-shrinker for Lagrangian mean curvature flow, and this is reflected in the fact that a circle in the plane will self-similarly shrink to a point under the flow \eqref{eq:curve.flow}.  Hence, the flow starting at the Clifford torus has a Type I singularity.
\end{example}

\begin{example}[Ellipses] Recently, it was shown in \cite{ELS}, that for any Lagrangian torus $L$ defined by an ellipse as in Figure \ref{fig:ellipse} with $a\neq b$, the Lagrangian mean curvature flow starting at $L$ must develop a finite-time singularity which is \emph{not} modelled on the Clifford torus.  Hence, the Clifford torus is \emph{unstable} under Lagrangian mean curvature flow, even under arbitrarily small Hamiltonian perturbations (which correspond to taking variations with $ab$ constant). 

Further, it is known by work in \cite{Grohetal,NevesMonotone} that, for $a\gg b$,   the Lagrangian mean curvature flow starting at a torus defined by such an ellipse would have to develop a Type II singularity at the origin, whose Type I blow-up is a circle-invariant Lagrangian pair of planes defined by a pair of lines in $\mathbb{C}$ which meet at right angles at the origin.  The curve at the singular time will look something like the figure eight curve in Figure \ref{fig:8} below, and thus the Lagrangian torus will become a 2-sphere at the singular time.  We expect that this same behaviour occurs for any $a\neq b$.
\end{example}

\begin{example}[Star-shaped curves]
In \cite{Grohetal}, the authors gave conditions under which a star-shaped curve with respect to the origin will contract to a point under the flow \eqref{eq:curve.flow} so that the corresponding Lagrangian mean curvature flow has a finite-time Type I singularity.  Moreover, the Type I blow-up is a Lagrangian self-shrinker constructed in \cite{Anciaux}.  However, in $\mathbb{C}^2$, the constraints on the star-shaped curve mean that the   Lagrangian tori cannot be embedded.
\end{example}

\begin{example}[Chekanov torus]
It is natural to ask about the behaviour of Lagrangian mean curvature flow for embedded tori defined by embedded closed curves which \emph{do not} enclose the origin. Such tori include the Chekanov torus \cite{Chekanov}, which is important in symplectic topology, particularly as it is monotone and \emph{not} Hamiltonian isotopic to the Clifford torus: see Figure \ref{fig:Chekanov} for a curve corresponding to an equivariant Chekanov torus.   

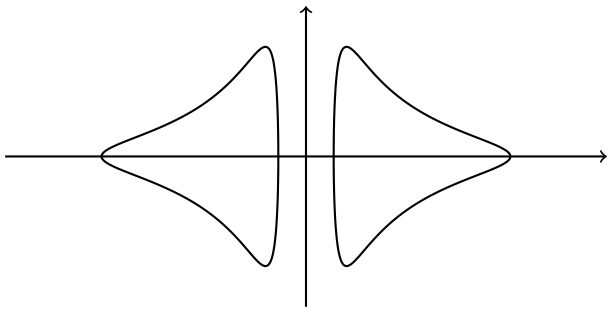
\begin{figure}[h]
\begin{center}
\begin{tikzpicture}
\draw[thick,->] (0,-2) -- (0,2);
\draw[thick,->] (-4,0) -- (4,0);
\draw[thick,domain=0:2*pi,samples=400] plot ({exp(cos(\x r)))},{sin(\x r)*exp(-cos(\x r))});
\draw[thick,domain=0:2*pi,samples=400] plot ({-exp(cos(\x r)))},{-sin(\x r)*exp(-cos(\x r))});
\end{tikzpicture}
\end{center}
\caption{A curve defining a Chekanov torus $T^2$ in $\mathbb{C}^2$}\label{fig:Chekanov}
\end{figure}

For some examples of Lagrangian tori defined by embedded closed curves not enclosing the origin, it is shown in \cite{Groh} that the Lagrangian mean curvature flow collapses to a circle in finite time, so the corresponding curve evolving under \eqref{eq:curve.flow} contracts to a point outside of the origin.  In these cases, the Lagrangian mean curvature flow has a Type I singularity whose Type I blow-up is a self-shrinking cylinder $S^1\times \mathbb{R}$.  However, the general picture is not know, and certainly seems worth exploring.
\end{example}

\section{Spheres: immersed closed curves}

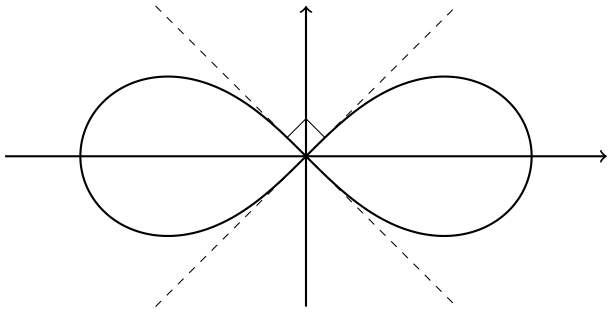
\begin{figure}[h]
\begin{center}
\begin{tikzpicture} 
\draw[thick,->] (0,-2) -- (0,2);
\draw[thick,->] (-4,0) -- (4,0);
\draw[dashed] (-2,-2) -- (2,2);
\draw[dashed] (-2,2) -- (2,-2);
\draw (-1/4,1/4) -- (0,1/2);
\draw (1/4,1/4) -- (0,1/2);
\draw[thick,domain=0:2*pi,samples=400] plot ({-3*sin(\x r)/(1+cos(\x r)*cos(\x r))}, {3*sin(\x r)*cos(\x r)/(1+cos(\x r)*cos(\x r))});
\end{tikzpicture}
\caption{A figure eight defining an immersed Lagrangian 2-sphere $S^2$}\label{fig:8}
\end{center}
\end{figure}

In this section we consider Lagrangian 2-spheres in $\mathbb{C}^2$ (which must be immersed) and are equivariant, so defined by curves as in Figure \ref{fig:8}.  

The curve $\gamma$ in Figure \ref{fig:8} and the resulting Lagrangian $L$ defined by $\gamma$ as in \eqref{eq:L.equiv} has two distinct properties.
\begin{enumerate}
\item[(a)] The curve $\gamma$ is contained in the region defined by the dashed lines, which each make an angle of $\frac{\pi}{4}$ with the horizontal axis, and $\gamma$ meets any circle centred at the origin in at most 4 points.
\item[(b)] The Ricci curvature of the induced metric on $L$ satisfies $\text{Ric}\geq cr^2$ where $c>0$ is a constant and $r$ is the distance to the origin in $\mathbb{C}^2$. 
\end{enumerate}
We shall describe results concerning curves satisfying each of these properties.

\begin{example}[Case (a)]
For curves $\gamma$ as in (a), Viana \cite{Viana} showed that the flow \eqref{eq:curve.flow} starting at $\gamma$ has a finite-time singularity at the origin where the flow shrinks to a point.  The corresponding Lagrangian mean curvature flow has a Type II singularity at the origin, whose Type I blow-up is a plane with multiplicity two (roughly, two copies of the same plane).
\end{example}

\begin{example}[Case (b)]
For curves $\gamma$ as in (b), it was shown in \cite{SmoczykWhitney} that again the flow \eqref{eq:curve.flow} starting at  $\gamma$ has a finite-time singularity at  $0$, where the flow shrinks to a point.  Moreover, the Lagrangian mean curvature flow given by $\gamma$ has a Type II singularity at the origin whose Type II blow-up is the product of a Grim Reaper curve and a real line.  The Grim Reaper curve is shown in Figure \ref{fig:grim} and is given by 
\begin{equation}\label{eq:grim}
\gamma(t)=\{-\log\cos y+iy\in\mathbb{C}\,:\,y\in(-\textstyle\frac{\pi}{2},\frac{\pi}{2})\}.
\end{equation}
It is a \emph{translator} for curve shortening flow, i.e.~it just translates to the right along the flow.
\begin{figure}[h]
\begin{center}
\begin{tikzpicture} [xscale=1,yscale=1]
\draw[thick,->] (0,-2) -- (0,2);
\draw[thick,->] (-1,0) -- (5,0);
\draw[dashed] (-1,pi/2) -- (5,pi/2);
\draw[dashed] (-1,-pi/2) -- (5,-pi/2);
\draw[thick,domain=-pi/2+0.006:pi/2-0.005,samples=400] plot ({-ln(cos(\x r))},{\x});
\end{tikzpicture}
\end{center}
\caption{Grim Reaper curve}\label{fig:grim}
\end{figure}
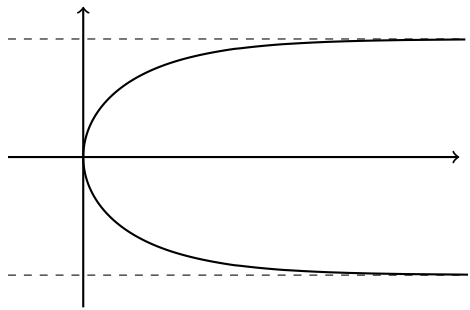
\end{example}

\begin{figure}[h]
\begin{center}
\begin{tikzpicture} 
\draw[thick,->] (0,-2) -- (0,2);
\draw[thick,->] (-4,0) -- (4,0);
\draw[dashed] (-2,-2) -- (2,2);
\draw[dashed] (-2,2) -- (2,-2);
\draw (-1/4,1/4) -- (0,1/2);
\draw (1/4,1/4) -- (0,1/2);
\draw[thick,domain=0:2*pi,samples=400] plot ({-2*sin(\x r)/(1+cos(\x r)*cos(\x r))}, {1*sin(\x r)*cos(\x r)/(1+cos(\x r)*cos(\x r))});
\end{tikzpicture}
\caption{A collapsing figure eight along the flow \eqref{eq:curve.flow}}\label{fig:8.collapse}
\end{center}
\end{figure}
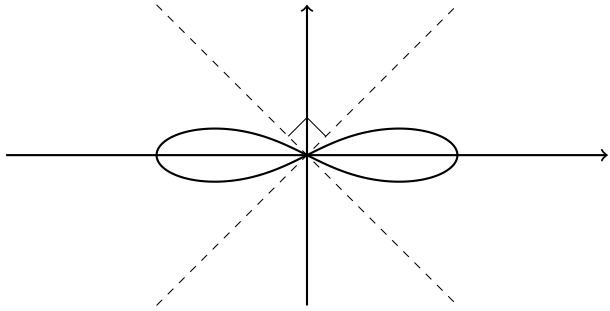

\begin{example}[Whitney sphere/figure eight]  
The figure eight  in Figure \ref{fig:8} defines the well-known Whitney sphere in $\mathbb{C}^2$.  The flow \eqref{eq:curve.flow} starting at this curve  cannot be self-similar as there are no Lagrangian self-shrinking spheres \cite{ChenMa}.  However, by the work in \cite{SmoczykWhitney,Viana} it will still shrink to a point in finite time.  

At the finite time singularity the Lagrangian mean curvature flow will have a Type II singularity whose Type I blow-up is a plane with multiplicity two, and the Type II blow-up is the product of a Grim Reaper curve with a line.  Along the flow it deforms as in Figure \ref{fig:8.collapse}, ``squashing'' vertically faster than it does horizontally, and one can imagine the Grim Reaper curve emerging at the extreme left (and right) points of the curve in the limit. 
\end{example}

\section{Cylinders: embedded open curves}

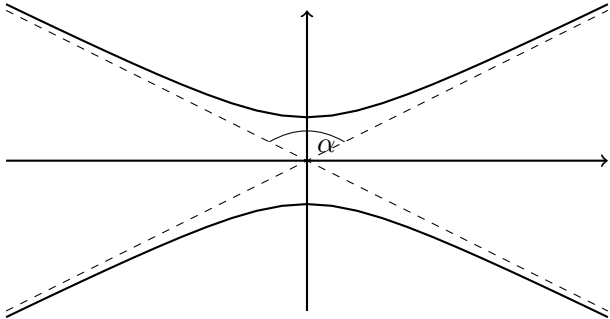
\begin{figure}[h]
\begin{center}
\begin{tikzpicture} 
\draw[thick,->] (0,-2) -- (0,2);
\draw[thick,->] (-4,0) -- (4,0);
\draw[dashed] (-4,-2) -- (4,2);
\draw[dashed] (-4,2) -- (4,-2);
\draw[thick,domain=-4:4]plot (\x,{sqrt(1/4*\x*\x+1/3)});
\draw[thick,domain=-4:4]plot (\x,{-sqrt(1/4*\x*\x+1/3)});
\draw (-2/4,1/4) to [out=30,in=150] (2/4,1/4);
\node[below right] at (0,5/12) {$\alpha$};
\end{tikzpicture}
\end{center}
\caption{A non-compact arc defining an asymptotically planar Lagrangian cylinder $S^1\times\mathbb{R}$}\label{fig:arc}
\end{figure}

For the final set of examples, we look at embedded non-compact curves asymptotic to straight lines, as in Figure \ref{fig:arc}.  The behaviour of the flow crucially depends on the angle $\alpha$ in Figure \ref{fig:arc}.

\begin{example}[$\alpha=\frac{\pi}{2}$: Lawlor necks]
In this case, the curve $\gamma$ in Figure \ref{fig:arc} is a standard hyperbola asymptotic to a pair of straight lines meeting at right angles.  The corresponding Lagrangian $L$ is minimal, i.e.~it has $H=0$, and so is a stationary point for Lagrangian mean curvature flow \eqref{eq:LMCF}.  The Lagrangian $L$ is called a Lawlor neck.

Recently, Su \cite{Su} considered a natural class of curves $\gamma$ asymptotic to a pair of lines with angle $\alpha=\frac{\pi}{2}$ as in Figure \ref{fig:arc}, which, in particular, are sandwiched between hyperbolae defining Lawlor necks.  Su showed that the Lagrangian mean curvature flow starting at the Lagrangian defined by $\gamma$ exist for all time and converges to a Lawlor neck.
\end{example}

\begin{example}[$\alpha\in(0,\frac{\pi}{2})$: self-expanders]
For every $\alpha\in(0,\frac{\pi}{2})$, Anciaux \cite{Anciaux} showed that there is a dilation family of curves $\gamma$ as in Figure \ref{fig:arc} so that the corresponding Lagrangian is a \emph{self-expander}; i.e.~it simply expands by dilation under the flow \eqref{eq:LMCF}.  

Su \cite{Su} considered certain natural curves $\gamma$ as in Figure \ref{fig:arc} for $\alpha\in(0,\frac{\pi}{2})$, bounded between two curves defining self-expanders, and showed that (after rescaling) Lagrangian mean curvature flow starting at the Lagrangian defined by $\gamma$ exists for all time and converges to an Anciaux self-expander.
\end{example}

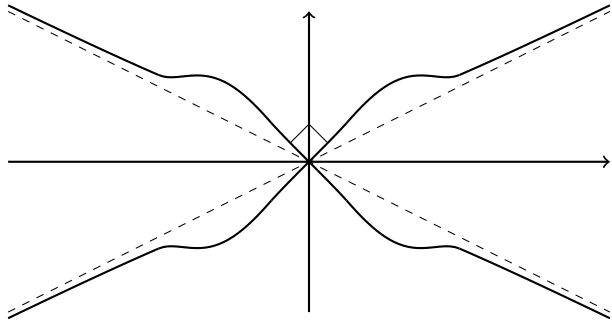
\begin{figure}[h]
\begin{center}
\begin{tikzpicture} 
\draw[thick,->] (0,-2) -- (0,2);
\draw[thick,->] (-4,0) -- (4,0);
\draw[dashed] (-4,-2) -- (4,2);
\draw[dashed] (-4,2) -- (4,-2);
\draw[thick,domain=-1/2:1/2]plot (\x,{abs(\x)*sqrt(1/4*\x*\x+1)});
\draw[thick,domain=1/2:2]plot (\x,{abs(\x)*sqrt(1/4*\x*\x+1)*16/15*16/15*(1-1/4*\x*\x)*(1-1/4*\x*\x)+(1-16/15*16/15*(1-1/4*\x*\x)*(1-1/4*\x*\x))*sqrt(1/4*\x*\x+1/3)});
\draw[thick,domain=-2:-1/2]plot (\x,{abs(\x)*sqrt(1/4*\x*\x+1)*16/15*16/15*(1-1/4*\x*\x)*(1-1/4*\x*\x)+(1-16/15*16/15*(1-1/4*\x*\x)*(1-1/4*\x*\x))*sqrt(1/4*\x*\x+1/3)});
\draw[thick,domain=-4:-2]plot (\x,{sqrt(1/4*\x*\x+1/3)});
\draw[thick,domain=2:4]plot (\x,{sqrt(1/4*\x*\x+1/3)});
\draw[thick,domain=-4:-2]plot (\x,{-sqrt(1/4*\x*\x+1/3)});
\draw[thick,domain=2:4]plot (\x,{-sqrt(1/4*\x*\x+1/3)});
\draw[thick,domain=-1/2:1/2]plot (\x,{-abs(\x)*sqrt(1/4*\x*\x+1)});
\draw[thick,domain=1/2:2]plot (\x,{-abs(\x)*sqrt(1/4*\x*\x+1)*16/15*16/15*(1-1/4*\x*\x)*(1-1/4*\x*\x)-(1-16/15*16/15*(1-1/4*\x*\x)*(1-1/4*\x*\x))*sqrt(1/4*\x*\x+1/3)});
\draw[thick,domain=-2:-1/2]plot (\x,{-abs(\x)*sqrt(1/4*\x*\x+1)*16/15*16/15*(1-1/4*\x*\x)*(1-1/4*\x*\x)-(1-16/15*16/15*(1-1/4*\x*\x)*(1-1/4*\x*\x))*sqrt(1/4*\x*\x+1/3)});
\draw (-1/4,1/4) -- (0,1/2);
\draw (1/4,1/4) -- (0,1/2);
\end{tikzpicture}
\end{center}
\caption{Singularity formation for an asymptotically planar Lagrangian cylinder}\label{fig:arc.2}
\end{figure}

\begin{example}[$\alpha\in(\frac{\pi}{2},\pi)$: singularities]
When $\alpha\in(\frac{\pi}{2},\pi)$, it is shown in \cite{NevesSingularities,Wood} that for curves $\gamma$ as in Figure \ref{fig:arc} the flow \eqref{eq:curve.flow} has a first finite-time singularity at the origin.  Moreover,  the singularity is Type II and the Type I blow-up is a circle-invariant pair of Lagrangian planes defined by a pair of lines in $\mathbb{C}$ which meet at right angles at $0$.  A rough picture of what happens at the finite-time singularity is given in Figure \ref{fig:arc.2}.  The Lagrangian has become two tranversely intersecting copies of $\mathbb{R}^2$ at the singular time.

Recently, Wood \cite{Wood} has shown further that the Type I blow-up at the origin is unique, the Type II blow-up is unique (up to scale) and given by the Lawlor neck asymptotic to the Type I blow-up, so that the blow-down of the Type II blow-up is equal to the Type I blow-up.  This matches well with the classification results for Type II blow-ups in \cite{LLS}.  
\end{example}

\section{Outlook}

There are many further questions one may ask about general Lagrangian mean curvature flow which may already yield an interesting answer in the equivariant setting, motivated in part by the examples we have discussed.  
\begin{itemize}
\item When can one relate blow-downs of Type II blow-ups to Type I blow-ups at a finite-time singularity?
\item Does the mean curvature necessarily blow-up at a finite-time singularity?
\item Can one better understand the role of pseudoholomorphic curves in the formation of singularities in Lagrangian mean curvature flow?
\end{itemize}
Answering any of these questions would be significant for further study of Lagrangian mean curvature flow.
 

\providecommand{\bysame}{\leavevmode\hbox to3em{\hrulefill}\thinspace}

\end{document}